\documentclass{amsart}

\newcommand{\R}{{\mathbb R}}
\newcommand{\N}{{\mathbb N}}
\newcommand{\on}{\operatorname}
\newcommand{\eps}{\varepsilon}
\newtheorem{theorem}{Theorem}[section]
\newtheorem{lemma}[theorem]{Lemma}
\newtheorem{corollary}[theorem]{Corollary}
\newtheorem{proposition}[theorem]{Proposition}

\theoremstyle{definition}

\theoremstyle{remark}
\newtheorem{remark}[theorem]{Remark}

\numberwithin{equation}{section}

\begin{document}

\title{Energy identity for approximations of harmonic maps from surfaces}

\author{Tobias Lamm}
\address{Max-Planck-Institute for Gravitational Physics, Am M\"uhlenberg 1, 14476 Golm, Germany}
\curraddr{Department of Mathematics, University of British Columbia, 1984 Mathematics Road, Vancouver, BC, V6T 1Z2, Canada}
\email{tlamm@math.ubc.ca}
\thanks{The author would like to thank Yuxiang Li for pointing out an error in an earlier version of the paper.}

%    General info

\date{December 12, 2007.}


\begin{abstract}
We prove the energy identity for min-max sequences of the Sacks-Uhlenbeck and the biharmonic approximation of harmonic maps from surfaces into general target manifolds. The proof relies on Hopf-differential type estimates for the two approximations and on estimates for the concentration radius of bubbles.
\end{abstract}

\maketitle

\section{Introduction}
Let $(M^2,g)$ be a smooth and compact Riemannian surface and let $(N^n,h)$ be a smooth and compact Riemannian manifold, both without boundary. We assume that $N^n\hookrightarrow \R^m$ isometrically. For $u\in W^{1,2}(M,N)$ we define the Dirichlet energy 
\begin{align}
E(u)=\int_M |\nabla u|^2 dv_g. \label{diri}
\end{align}
Critical points of $E$ are called harmonic maps and they solve the elliptic system
\begin{align}
\Delta u+A(u)(\nabla u, \nabla u)=0, \label{harm}
\end{align}
where $A$ is the second fundamental form of the embedding $N\hookrightarrow \R^m$. The geometric interest in harmonic maps from surfaces comes from the fact that if the harmonic map is additionally conformal (i.e. angle preserving) then the image of the map is a minimal immersion of $M$ in $N$. For example it is well known that every harmonic map $u:S^2\rightarrow N$ is minimal. It is therefore of interest to find critical points of the Dirichlet energy. Since $E$ does not satisfy the Palais-Smale condition the classical variational methods do not apply to $E$. In order to overcome this difficulty Sacks \& Uhlenbeck \cite{sacks81} introduced a regularization of the Dirichlet energy. More precisely, they considered for every $\alpha>1$ and $u\in W^{1,2\alpha}(M,N)$ the functional
\begin{align}
E_\alpha(u)=\int_M (1+|\nabla u|^2)^\alpha dv_g. \label{SUenergy}
\end{align}
Since this functional satisfies the Palais-Smale condition they were able to show the existence of a smooth critical point of $E_\alpha$ for every $\alpha>1$ by classical variational methods. These critical points $u_\alpha$ solve the elliptic system
\begin{align}   
\frac{1}{\alpha}dE_\alpha(u)=\on{div}((1+|\nabla u_\alpha|^2)^{\alpha-1}\nabla u_\alpha)+(1+|\nabla u_\alpha|^2)^{\alpha-1}A(u_\alpha)(\nabla u_\alpha,\nabla u_\alpha)=0.\label{SUeq}
\end{align}
Sacks \& Uhlenbeck then studied sequences of critical points $u_\alpha$ ($\alpha \rightarrow 1$) of $E_\alpha$ with uniformly bounded energy $E_\alpha(u_\alpha)\le c$. They showed that for a subsequence $\alpha_k\rightarrow 1$ the maps $u_{\alpha_k}$ converge weakly in $W^{1,2}(M,N)$ and strongly away from at most finitely many singular points to a smooth harmonic map $u_1 \in C^\infty(M,N)$. Moreover they were able to perform a blow-up around these finitely many singular points and they showed that the blow-up`s are non-trivial minimal two-spheres. As an application of this analysis Sacks \& Uhlenbeck proved the existence of a minimal two-sphere in every homotopy class if $\pi_2(N)=0$. 

What was left over in their analysis of sequences of critical points of $E_\alpha$ was the question if there is some energy-loss occurring during the blow-up process.

In \cite{lamm06} the author considered a different regularization of the Dirichlet energy, namely for every $\eps>0$ and every $u\in W^{2,2}(M,N)$ we studied the functional 
\begin{align}     
E_\eps (u) =\int_M |\nabla u|^2 dv_g +\eps \int_M |\Delta u|^2 dv_g. \label{fourth1}
\end{align}
The Euler-Lagrange equation of $E_\varepsilon$ is given by 
\begin{align}
\Delta u-\varepsilon \Delta^2 u=-A(u)(\nabla u,\nabla u)+f[u], \label{fourth2}
\end{align}
where $f[u]\perp T_uN$ and
\begin{align}
|f[u]| \le c(|u|) \varepsilon (|\nabla u| \nabla^3 u|+|\nabla^2 u|^2+|\nabla u|^4). \label{fourth3}
\end{align}
For every $\varepsilon>0$ the functional $E_\varepsilon$ satisfies the Palais-Smale condition and therefore critical points exist and they are smooth. Hence, as in the case of the Sacks-Uhlenbeck approximation, we studied sequences $u_\eps \in C^\infty (M,N)$ ($\eps \rightarrow 0$) of critical points of $E_\eps$ with uniformly bounded energy $E_\eps(u_\eps)\le c$. We were able to show that for a subsequence $\eps_k \rightarrow 0$ the maps $u_{\eps_k}$ converge weakly in $W^{1,2}(M,N)$ and strongly away from at most finitely many singular points to a smooth harmonic map $u_0:M\rightarrow N$. Moreover, by performing a blow-up around the singular points, we showed that at most finitely many minimal two-spheres were separating. Additionally we were able to show that there is no energy lost during the blow-up process if $N=S^n\hookrightarrow \R^{n+1}$. The case of a general target manifold was left open.

In the main result of this paper we show that for both approximations and general target manifolds there is no energy-loss occurring if we assume an additional entropy-type condition. More precisely we have the following 
\begin{theorem}\label{main}
Let $(M^2,g)$ be a smooth, compact Riemannian surface without boundary and let $N$ be a smooth and compact Riemannian manifold without boundary, which we assume to be isometrically embedded into $\R^n$. Let $u_\alpha\in C^\infty(M,N)$ ($\alpha\to 1$) be a sequence of critical points of $E_\alpha$ with uniformly bounded energy. Moreover we assume that $u_\alpha$ satisfies 
\begin{align}
\liminf_{\alpha \rightarrow 1} (\alpha -1) \int_M \log (1+|\nabla u_\alpha|^2) (1+|\nabla u_\alpha|^2)^\alpha dv_g=0.\label{ent}
\end{align}
Then there exists a sequence $\alpha_k \to 1$ and at most finitely many points $x^1,\ldots, x^l \in M$ such that $u_{\alpha_k}\to u_1$ weakly in $W^{1,2}(M,N)$ and in $C^\infty_{\on{loc}}(M\backslash \{x^1,\ldots,x^l\},N)$ where $u_1:M\to N$ is a smooth harmonic map.

By performing a blow-up at each $x^i$, $1\leq i\leq l$, one gets that there exist at most finitely many non-trivial smooth harmonic maps $\omega^{i,j}:S^2\to N$, $1\leq j\leq j_i$, sequences of points $x_k^{i,j} \in M$, $x_k^{i,j} \to x^i$, and sequences of radii $r_k^{i,j}\in \R_+$, $r_k^{i,j} \to 0$, such that
\begin{align}
 \on{max} \{  \frac{r_k^{i,j}}{r_k^{i,j'}},\frac{r_k^{i,j'}}{r_k^{i,j}},\frac{\on{dist}(x_k^{i,j},x_k^{i,j'})}{r_k^{i,j}+r_k^{i,j'}} \} &\to \infty, \ \ \ \forall \ \ 1\leq i \leq l, \ \ 1\leq j,j'\leq j_i,\ \ j\not= j', \label{conv}\\   
\limsup_{k\rightarrow \infty} (r^{i,j}_k)^{1-\alpha_k} &=1 \ \ \ \forall \ \ 1\leq i\leq l,\ \ 1\leq j\leq j_i\ \ \ \text{and} \label{epsnullgeneral}\\
\lim_{k\to \infty}  E_{\alpha_k}(u_{\alpha_k}) &=E(u_1)+\on{vol}(M)+\sum_{i=1}^l\sum_{j=1}^{j_i} E(\omega^{i,j}). \label{energyequality}
\end{align}
\end{theorem} 
\begin{remark}
The Theorem remains true if we replace everywhere $E_\alpha$ by $E_\eps$, $u_\alpha$ by $u_\eps$, $u_1$ by $u_0$, the assumption \eqref{ent} by
\begin{align}
\liminf_{\eps \rightarrow 0} \eps \log(\frac{1}{\eps}) \int_M |\Delta u_\eps|^2 dv_g =0, \label{epsent}
\end{align}
the estimate \eqref{epsnullgeneral} by
\begin{align}
\limsup_{k\rightarrow \infty} \frac{\eps_k}{(r^{i,j}_k)^2} =0 \ \ \ \forall \ \ 1\leq i\leq l,\ \ 1\leq j\leq j_i, \label{epsnullgenerala}
\end{align}
and \eqref{energyequality} by
\begin{align}
\lim_{k\to \infty} E_{\eps_k}(u_{\eps_k})&=E(u_0)+\sum_{i=1}^l\sum_{j=1}^{j_i} E(\omega^{i,j}). \label{energyequalitya}
\end{align}
\end{remark}
\begin{remark}
By the results of Duzaar \& Kuwert \cite{duzaar98} (Theorem $2$) the above Theorem implies that we also have a decomposition in terms of homotopy classes.
\end{remark} 
Of course now one has to ask if there exist sequences of critical points of $E_\alpha$, resp. $E_\eps$, satisfying \eqref{ent}, resp. \eqref{epsent}. The answer to this question is yes and more precisely we have the following
\begin{lemma}\label{entropy}
Let $\alpha>1$ and let $\mathcal{F} \subset \mathcal{P}(W^{1,2\alpha}(M,N))$ be a collection of sets. Let $\Phi : [0,\infty[\times W^{1,2\alpha}(M,N) \rightarrow W^{1,2\alpha}(M,N)$ be any continuous semi-flow such that $\Phi(0,\cdot)=id$, $\Phi(t,\cdot)$ is a homeomorphism of $W^{1,2\alpha}(M,N)$ for any $t\ge 0$ and $E_\alpha(\Phi(t,u))$ is nonincreasing in $t$ for any $u\in W^{1,2\alpha}(M,N)$. We assume that $\Phi(t,F)\subset F$ for all $t\in [0,\infty)$ and all $F\in \mathcal{F}$. We define
\begin{align}
\beta_\alpha=\inf_{F\in \mathcal{F}} \sup_{u\in F} E_\alpha(u) \label{minmax}
\end{align}
and we assume that $\beta_\alpha <\infty$. Then for almost every $\alpha$ there exists a critical point $u_\alpha \in C^\infty(M,N)$ of $E_\alpha$ with $E_\alpha(u_\alpha)=\beta_\alpha$ and such that
\begin{align}
\liminf_{\alpha \rightarrow 1} (\alpha-1)\log(\frac{1}{\alpha-1})\int_M \log(1+|\nabla u_\alpha|^2)(1+|\nabla u_\alpha|^2)^\alpha dv_g =0.\label{mini}
\end{align}
With the obvious modifications the same conclusion remains true for the energy $E_\eps$. 
\end{lemma} 
\begin{remark}\label{examples}
For examples of subsets $\mathcal{F} \subset \mathcal{P}(W^{1,2\alpha}(M,N))$ satisfying the hypothesis of the above Lemma we refer the reader to \cite{palais70} (p. 190) or \cite{struwe00} (p. 88).
\end{remark}
As a Corollary of the above Theorem and Lemma, we obtain a new proof of a result of Jost \cite{jost91} on the energy identity for min-max sequences for the Dirichlet energy. 
\begin{corollary}\label{coro}
Let $(M^2,g)$ be a smooth, compact Riemannian surface without boundary and let $N\hookrightarrow \R^n$ be a smooth and compact Riemannian manifold without boundary. Moreover let $A$ be a compact parameter manifold, for simplicity we assume $\partial A=\emptyset$, and let $h_0:M\times A \to N$ be continuous. Let $H$ be the class of all maps homotopic to $h_0$ and
\begin{align}
 \beta:= \inf_{h\in H} \sup_{t\in A} E(h(\cdot,t)). \label{defbeta}
\end{align}
Then there exists a sequence $u_{\alpha_k} \in C^\infty (M,N)$ of critical points of $E_{\alpha_k}$, a harmonic map $u_1:M\to N$ and at most finitely many points $x^1,\ldots,x^l\in M$ such that 
\begin{align}
E_{\alpha_k}(u_{\alpha_k})=\beta_{\alpha_k}&= \inf_{h\in H} \sup_{t\in A} E_{\alpha_k}(h(\cdot,t)), \label{critical}\\
\beta_{\alpha_k} &\to \beta+\on{vol}(M), \label{convbeta}\\
u_{\alpha_k} &\rightharpoonup u_1\ \ \on{weakly}\ \ \on{in}\ \ W^{1,2}(M, N)\ \ \on{and} \label{weakkonvergenz}\\
u_{\alpha_k} &\to u_1 \ \ \on{in}\ \ C^\infty_{\on{loc}}(M\backslash \{x^1,\ldots,x^l\},N). \label{starkekonvergenz}
\end{align}
Moreover there exist at most finitely many non-trivial smooth harmonic maps $\omega^{i,j}:S^2\to N$, $1\leq i \leq l$, $1\leq j\leq j_i$, sequences of points $x_k^{i,j} \in M$, $x_k^{i,j} \to x^i$, and sequences of radii $r_k^{i,j}\in \R_+$, $r_k^{i,j} \to 0$, such that
\begin{align} 
\on{max} \{  \frac{r_k^{i,j}}{r_k^{i,j'}},\frac{r_k^{i,j'}}{r_k^{i,j}},\frac{\on{dist}(x_k^{i,j},x_k^{i,j'})}{r_k^{i,j}+r_k^{i,j'}} \} &\to \infty, \ \ \ \forall \ \ 1\leq i \leq l, \ \ 1\leq j,j'\leq j_i,\ \ j\not= j',
 \label{conv1a}\\   
\limsup_{k\rightarrow \infty} (r^{i,j}_k)^{1-\alpha_k} &=1\ \ \ \forall \ \ 1\leq i\leq l,\ \ 1\leq j\leq j_i\ \ \ \text{and} \label{epsnullgeneral1}\\
\lim_{k\to \infty} E_{\alpha_k}(u_{\alpha_k})&=E(u_1)+\on{vol}(M)+\sum_{i=1}^l\sum_{j=1}^{j_i} E(\omega^{i,j}). \label{energyequality10}
\end{align}
\end{corollary}
\begin{remark} 
With the obvious modifications the corollary remains true for the biharmonic approximation $E_\eps$.  
\end{remark}
\begin{proof}
The proof of this result is quite standard but we include it here for sake of completeness. It is obvious that for all $\alpha>1$ we have
\begin{align*}
 \beta +\on{vol}(M) \leq \beta_\alpha.
\end{align*}
Let $\delta >0$ and choose $\tilde{h} \in H \cap C^\infty(M\times A,N)$ such that
\begin{align*}
 \sup_{t\in A} E (\tilde{h}(\cdot,t)) \leq \beta+\delta .
\end{align*}
Then for $(\alpha-1)$ small enough we have
\begin{align*}
\sup_{t\in A} E_\alpha (\tilde{h}(\cdot,t)) &\leq \beta+\on{vol}(M)+\delta +c(\tilde{h})(\alpha-1)\\
 &\leq \beta +\on{vol}(M)+2\delta.
\end{align*}
This implies  
\begin{align*}
 \on{lim}_{\alpha \to 1} \beta_\alpha=\beta+\on{vol}(M).
\end{align*}
The result now follows from the minimax principle (see \cite{struwe00}), Theorem \ref{main} and Lemma \ref{entropy}.
\end{proof}
In the existing literature there are already some partial results available for the energy identity for the Sacks-Uhlenbeck approximation and there are many more results available for related problems. In the following we want to mention some of these results.

For the Sacks-Uhlenbeck approximation Duzaar \& Kuwert \cite{duzaar98} and Chen \& Tian \cite{chen99} proved the energy identity for sequences of minimizers of the energy $E_\alpha$ in a given homotopy class. Recently Moore \cite{moore06} proved the energy identity (he actually proved \eqref{energyequality} with the Dirichlet energy $E$ instead of the full $\alpha$-energy $E_\alpha$ on the left hand side) for min-max sequences of the Sacks-Uhlenbeck approximation under the additional assumption that the target manifold has finite fundamental group. The additional assumptions made by Chen \& Tian and Moore were used to ensure that the sequence of minimizers, respectively min-max sequence, converges to a geodesic of finite length on the necks connecting the bubbles and the weak limit (or body map) which then implies the energy identity. In our proof we use completely different arguments but we want to mention that it is not directly clear from our analysis that the sequence of critical points satisfying the entropy condition converges to a geodesic of finite length on the necks.

In a recent independent work, Li \& Wang \cite{li07} proved Theorem \ref{main} in the special case of sequences of minimizers (in their own homotopy class) of $E_\alpha$.

For sequences of harmonic maps and min-max sequences for the Dirichlet energy the energy identity was proved by Jost \cite{jost91} (see also \cite{parker96} for an alternative proof of the energy identity for sequences of harmonic maps).

Recently Colding \& Minicozzi \cite{colding07} proved the energy identity for sequences of maps with bounded Dirichlet energy which are "almost" conformal and which satisfy a certain replacement property.

The energy identity for the harmonic map heat flow and Palais-Smale sequences for the Dirichlet energy with tension field bounded in $L^2$ was established by Qing \cite{qing95} (in the case $N=S^n$) and independently by Ding \& Tian \cite{ding95} and Wang \cite{wang96} in the general case. Alternative proofs have been given by Qing \& Tian \cite{qing97} and Lin \& Wang \cite{lin98}. See also the paper of Topping \cite{topping04} for more refined results in this case.  

Lin \& Wang \cite{lin99}, \cite{lin02} used a Ginzburg-Landau approximation to regularize the Dirichlet energy and proved the energy identity in this situation. The disadvantage of the Ginzburg-Landau approximation is that the approximating maps do not have to map into the target manifold, only in the limit they are forced to do this. 

For maps from higher dimensional domains the energy identity for sequences of harmonic maps has been proved by Lin \& Rivi{\`e}re \cite{lin02b} for $N=S^n$. For other related problems such as sequences of Yang-Mills fields on a four-dimensional manifold, respectively biharmonic maps from a four-dimensional manifold into the sphere the energy identity has been proved by Rivi{\`e}re \cite{riviere02}, respectively Wang \cite{wang3}.

In the following we give a brief outline of the paper. 

In section $2$ we prove Theorem \ref{main} for the Sacks-Uhlenbeck approximation of harmonic maps. We start by recalling the small-energy regularity estimates and the blow-up procedure of Sacks \& Uhlenbeck \cite{sacks81} in section $2.1$. In Proposition \ref{concent} we prove the very important estimate for the concentration radius of the bubbles. The advantages of having a good estimate for the concentration radius can also be seen in the paper of Topping \cite{topping04}. In the next two sections we prove a Hopf-differential type estimate and an estimate for the tangential component of solutions of \eqref{SUeq} on annular regions. These estimates are proved in the same way as the corresponding estimates for harmonic maps, see for example \cite{sacks81} and \cite{ding95}. In section $2.4$ we use the bubbling induction argument of Ding \& Tian \cite{ding95} to reduce the proof of the energy identity to the case of one bubble. In this situation we then combine the previous estimates with the estimate for the concentration radius to complete the proof of the energy identity.

In section $3$ we treat the case of the biharmonic approximation. For this approximation the estimate for the concentration radius (see \eqref{concbih}) has already been proved in \cite{lamm06}. In section $3.1$ we review the small-energy estimates and the blow-up process from \cite{lamm06}. In section $3.2$ we use the stress-energy tensor of $E_\eps$ to get a Hopf-differential type estimate for the biharmonic approximation. The rest of the proof of the energy identity then follows as in the case of the Sacks-Uhlenbeck approximation and in the sections $3.3$ and $3.4$ we briefly describe the necessary modifications.

In section $4$ we use variational methods to prove Lemma \ref{entropy}. We follow closely the work of Struwe \cite{struwe00a}.

We use the notation $o_k(1)$, $o_{R_0}(1)$ and $o_R(1)$ to denote terms which tend to zero as $k\rightarrow \infty$, $R_0\rightarrow 0$ and $R\rightarrow \infty$ respectively.
\section{Energy identity for the Sacks-Uhlenbeck approximation of harmonic maps}

In this section we prove Theorem \ref{main} for the Sacks-Uhlenbeck approximation of harmonic maps.
\subsection{Results of Sacks and Uhlenbeck and estimates for the concentration radius}
We consider sequences of critical points $u_\alpha\in C^\infty (M,N)$ of the functional $E_\alpha$ with uniformly bounded energy $E_\alpha(u_\alpha)\le c$ and which satisfy the condition \eqref{ent}.  Due to the uniform boundedness of the energy it is easy to see that there exists a subsequence $\alpha_k\rightarrow 1$ such that 
\begin{align}
(\alpha_k -1) \int_M \log (1+|\nabla u_{\alpha_k}|^2) (1+|\nabla u_{\alpha_k}|^2)^{\alpha_k} dv_g\rightarrow 0.\label{enta}
\end{align}
and $u_{\alpha_k} \rightharpoonup u_1$ weakly in $W^{1,2}$. In section $3$ of \cite{sacks81} Sacks \& Uhlenbeck proved the following small energy regularity result for solutions of \eqref{SUeq}.
\begin{theorem}\label{regSU}
There exists $\varepsilon_0>0$ such that if $u_\alpha$ ($\alpha$ close to one) is a critical point of $E_\alpha$ with $\int_{B_{2R}}|\nabla u_\alpha|^2<\varepsilon_0$ (where $R>0$) then we have for every $m\in \N$
\begin{align}
\on{osc}_{B_R} u_\alpha+||\nabla^m u_\alpha||_{L^\infty(B_R)}R^m \le c(\int_{B_{2R}}|\nabla u_\alpha|^2)^{\frac{1}{2}}.\label{regSU1}
\end{align}
\end{theorem}
With the help of this Theorem Sacks \& Uhlenbeck were able to show that the sequence $u_{\alpha_k}$ converges strongly to a smooth harmonic map $u_1:M\rightarrow N$ away from finitely many points. These finitely many singular points $x^i\in M$, $1\le i \le l$, are caracterized by the condition that
\begin{align}
\limsup_{k\rightarrow \infty} E(u_{\alpha_k},B_R(x^i))\ge \varepsilon_0,\label{SUconc}
\end{align}
for every $R>0$ and every $1\le i\le l$. Around these finitely many singular points they were able to perform a blow-up and show that a non-trivial harmonic two-sphere separates. The blow-up can be done as follows:  Fix $R_0>0$ such that $B_{R_0}(x^i)\cap B_{R_0}(x^j)=\emptyset$ for every $i,j\in \{1,\ldots,l\}$, $i\not=j$. Because of \eqref{SUconc} there exists a sequence of points $x^i_k\rightarrow x^i$ and radii $r_k^i\rightarrow 0$ such that 
\begin{align}
\max_{y\in B_{R_0}(x^i)}E(u_{\alpha_k},B_{r_k^i}(y))=E(u_{\alpha_k},B_{r_k^i}(x^i_k))=\frac{\varepsilon_0}{2}.\label{SUcfct}
\end{align}
Defining:
\begin{align}
v^i_k:&B_{\frac{R_0}{r_k^i}}\rightarrow N \nonumber \\
v^i_k(x)=& u_{\alpha_k}(x^i_k+r_k^ix) \label{res}
\end{align}
we see that $v^i_k$ solves \eqref{SUeq} with $1$ replaced by $(r_k^i)^2$ and moreover
\begin{align}
\max_{y\in B_{\frac{R_0}{2r_k^i}}}E(v^i_k,B_{1}(y))=E(v^i_k,B_{1}(0))=\frac{\varepsilon_0}{2}.\label{SUcfct1}
\end{align}
Therefore we can apply Theorem \ref{regSU} to $v^i_k$ and get that $v^i_k$ converges in $C^1$ to a smooth harmonic map $\omega^i$ from $\R^2$ into $N$. By the point removabilty result of Sacks \& Uhlenbeck we can then extend $\omega^i$ to a smooth harmonic map from $S^2$ to $N$. 
\newline 
As a consequence of this blow-up procedure we get the following estimate for the concentration radius .
\begin{lemma}\label{concentrationradius}
Using the above notation we have that
\begin{align}
1\le \limsup_{k\rightarrow \infty} (r^i_k)^{1-\alpha_k}<\infty,\label{concentrationradius1}
\end{align}
for every $1\le i \le l$.
\end{lemma}
\begin{proof}
Because of \eqref{SUcfct} and H\"older's inequality we know that
\begin{align*}
\frac{\varepsilon_0}{2}&=E(u_{\alpha_k},B_{r^i_k}(x^i_k))\\
&\le (\int_M (1+|\nabla u_{\alpha_k}|^2)^{\alpha_k})^{\frac{1}{\alpha_k}}(r^i_k)^{\frac{2(\alpha_k-1)}{\alpha_k}}\\
&\le c(r^i_k)^{\frac{2(\alpha_k-1)}{\alpha_k}}.
\end{align*}
From this the claim follows.
\end{proof}
In the next Proposition we use \eqref{enta} to improve the above estimate for the concentration radius (see also \cite{struwe00a} were this was observed for a similar approximation of a different problem).
\begin{proposition}\label{concent}
We have that
\begin{align}
\lim_{k\rightarrow \infty} (r_k^i)^{1-\alpha_k}=1,\label{concentrationradiusent}
\end{align}
for every $1\le i \le l$.
\end{proposition}
\begin{proof}
We let $\varepsilon_0$ be as above and we assume without loss of generality that $l=1$. Furthermore we let $r_k^1=r_k$, $x_k^1=x_k$ and $u_{\alpha_k}=u_k$. For every $k\in \N$ we define the set
\begin{align}
\Omega_k =\{ x\in B_{r_k}(x_k)| |\nabla u_k(x)|\ge \frac{\sqrt{\eps_0}}{2\sqrt{\pi}r_k}  \} \label{concent1}
\end{align}
and we claim that there exists a constant $c>0$ such that for every $k\in \N$ we have
\begin{align}
|\Omega_k| \ge cr_k^2. \label{concent2}
\end{align}
If this is not the case we can find a subsequence $k_m$ such that
\begin{align}
|\Omega_{k_m}| \le \frac{r_{k_m}^2}{m}. \label{concent3}
\end{align}
From \eqref{SUcfct} and Theorem \ref{regSU} we get 
\begin{align}
||\nabla u_{k_m}||_{L^\infty(B_{r_{k_m}}(x_{k_m}))} \le \frac{c\sqrt{\varepsilon_0}}{r_{k_m}}. \label{estv}
\end{align}
From the definition of $\Omega_{k_m}$ we see that for every $x\in B_{r_{k_m}}(x_{k_m})\backslash \Omega_{k_m}$ we have the estimate
\begin{align}
|\nabla u_{k_m}|(x)\le \frac{\sqrt{\eps_0}}{2\sqrt{\pi}r_{k_m}}.\label{estv2}
\end{align}
Using \eqref{concent3}, \eqref{estv} and \eqref{estv2} we get from \eqref{SUcfct}
\begin{align*}
\frac{\varepsilon_0}{2} =& \int_{B_{r_{k_m}}(x_{k_m})} |\nabla u_{k_m}|^2  = \int_{\Omega_{k_m}}|\nabla u_{k_m}|^2+\int_{B_{r_{k_m}}(x_{k_m})\backslash \Omega_{k_m}}|\nabla u_{k_m}|^2\\
\le& \frac{c\varepsilon_0 |\Omega_{k_m}|}{r_{k_m}^2}+\pi r_{k_m}^2 \frac{\eps_0}{4\pi r_{k_m}^2}\\
\le& \frac{c}{m}+\frac{\eps_0}{4}\\
\rightarrow& \frac{\eps_0}{4} ,
\end{align*}
as $m\rightarrow \infty$. This contradiction proves the estimate \eqref{concent2}.

Now we use \eqref{enta}, Lemma \ref{concentrationradius}, the definition of $\Omega_k$ and \eqref{concent2} to estimate 
\begin{align*}
0=& \lim_{k\rightarrow \infty} (\alpha_k-1)\int_M \log (1+|\nabla u_k|^2)(1+|\nabla u_k|^2)^{\alpha_k}\\
\ge& \lim_{k\rightarrow \infty} (\alpha_k-1)\int_{\Omega_k} \log (|\nabla u_k|^2)|\nabla u_k|^{2\alpha_k}\\  
\ge&  c\lim_{k\rightarrow \infty} (\alpha_k-1) r_k^{2(1-\alpha_k)}\log(\frac{\eps_0}{4\pi r_k^2})\\
=& c\lim_{k\rightarrow \infty} (1-\alpha_k) r_k^{2(1-\alpha_k)}(2\log r_k-\log \eps_0+\log 4\pi)\\
=& c\lim_{k\rightarrow \infty} r_k^{2-2\alpha_k} \log r_k^{2-2\alpha_k}\\
\ge& c\lim_{k\rightarrow \infty}\log r_k^{2-2\alpha_k}\\
\ge& 0
\end{align*}
and hence the desired convergence result for the concentration radius follows.
\end{proof}
\subsection{A Hopf differential type estimate}
In the case of sequences of harmonic maps or Palais-Smale sequences for the Dirichlet energy with tension field bounded in $L^2$ an important ingredient in the proof of the energy identity was an estimate for the Hopf differential (see e.g. \cite{ding95}, \cite{sacks81}). In the next lemma we show that a related result is true for solutions of \eqref{SUeq}.
\begin{lemma}\label{radialSU}
Let $B\subset \R^2$ be the unit ball and let $u_\alpha\in C^\infty(B,N)$ be a solution of \eqref{SUeq}. Then we have for every $0<r<1$ and every $\alpha$ close to one
\begin{align}
\int_{\partial B_r} (1+|\nabla u_\alpha|^2)^{\alpha-1} |(u_\alpha)_r|^2\le& c\int_{\partial B_r}(1+\frac{|(u_\alpha)_\theta|^2}{r^2}) (1+|\nabla u_\alpha|^2)^{\alpha-1} \nonumber \\
&+\frac{c(\alpha-1)}{r}\int_{B_r} (1+|\nabla u_\alpha|^2)^{\alpha}. \label{radialSU1}
\end{align}
\end{lemma}
\begin{proof}
We multiply equation \eqref{SUeq} by $x\cdot \nabla u_\alpha$ and integrate over $B_r$ to get (remember that $A(u_\alpha)(\nabla u_\alpha,\nabla u_\alpha) \perp T_{u_\alpha} N$ for every $x\in B$)
\begin{align*}
0=& \int_{B_r} \on{div}((1+|\nabla u_\alpha|^2)^{\alpha-1}\nabla u_\alpha)x\cdot \nabla u_\alpha\\
=& -\int_{B_r}(1+|\nabla u_\alpha|^2)^{\alpha-1}|\nabla u_\alpha|^2+\int_{\partial B_r}r(1+|\nabla u_\alpha|^2)^{\alpha-1}|(u_\alpha)_r|^2\\
&-\frac{1}{2}\int_{B_r}(1+|\nabla u_\alpha|^2)^{\alpha-1}x\cdot \nabla (1+|\nabla u_\alpha|^2).
\end{align*}
Next we integrate by parts and get
\begin{align*}
\frac{\alpha}{2}\int_{B_r}(1+|\nabla u_\alpha|^2)^{\alpha-1}x\cdot \nabla (1+|\nabla u_\alpha|^2)=& -\int_{B_r}(1+|\nabla u_\alpha|^2)^{\alpha}\\
&+\int_{\partial B_r}\frac{r}{2}(1+|\nabla u_\alpha|^2)^{\alpha}.
\end{align*}
Using the identity 
\begin{align*}
|\nabla u_\alpha|^2=|(u_\alpha)_r|^2+\frac{1}{r^2}|(u_\alpha)_\theta|^2
\end{align*}
and combining everything we end up with
\begin{align*}
\int_{\partial B_r}r(1+|\nabla u_\alpha|^2)^{\alpha-1}|(u_\alpha)_r|^2\le& c\int_{\partial B_r}\frac{1}{r}(1+|\nabla u_\alpha|^2)^{\alpha-1}|(u_\alpha)_\theta|^2\\
&+c(\alpha-1)\int_{B_r}(1+|\nabla u_\alpha|^2)^{\alpha}\\
&+c\int_{\partial B_r}r(1+|\nabla u_\alpha|^2)^{\alpha-1}.
\end{align*}
\end{proof}
\subsection{Estimate for the tangential component}
In this section we show that if the Dirichlet energy is small on all annular regions with bounded geometry then the tangential derivative of $u_\alpha$ converges to zero on the annular region which is the union of all the annuli with bounded geometry. The proof of this fact follows closely the previous work of Sacks \& Uhlenbeck \cite{sacks81} and Ding \& Tian \cite{ding95}. In the following we use for $0<a_1<a_2<1$ the notation $A(a_1,a_2)=\{ x\in \R^2 |a_1 \le |x| \le a_2\}$.
\begin{lemma}\label{SUtangential}
There exists $\delta_0>0$ such that for all $\delta < \delta_0$ and all solutions $u_\alpha\in C^\infty(B,N)$ of \eqref{SUeq} with $\int_{A(r,2r)}|\nabla u_\alpha|^2 < \delta$ for every $r\in (R_1,\frac{R_2}{2})$, we have for $\alpha-1$ small enough
\begin{align}
\int_{2R_1}^{\frac{R_2}{4}}\int_0^{2\pi} \frac{1}{r} |(u_\alpha)_\theta|^2dr d\theta \le c\sqrt{\delta}(1+(\log R_1^{1-\alpha})). \label{SUtangential1}
\end{align}
\end{lemma}
\begin{proof}
Let $\delta_0<\varepsilon_0$ and let $y\in A(2R_1,\frac{R_2}{4})$. Then we have that $\frac{2|y|}{3},\frac{4|y|}{3} \in (R_1,\frac{R_2}{2})$ and $B_{\frac{|y|}{3}}(y)\subset B_{\frac{4|y|}{3}}\backslash B_{\frac{2|y|}{3}}$. From our assumption and Theorem \ref{regSU} we therefore conclude that 
\begin{align}
\sum_{i=1}^2 |x|^i |\nabla^i u_\alpha|(x) \le c\sqrt{\delta}, \label{SUtangential2}
\end{align}
for every $x\in A(2R_1,\frac{R_2}{4})$. Now we let $\frac{R_2}{4R_1}=2^l+q$, $l\in \N$ and $q\ge 0$, and define $A_k=A(2^k R_1, 2^{k+1}R_1)$ for all $1\le k\le l-1$ and we let $A_{l}=A(2^l R_1,\frac{R_2}{4})$. Next we note that equation \eqref{SUeq} can equivalently be written as
\begin{align}
\Delta u_\alpha+A(u_\alpha)(\nabla u_\alpha, \nabla u_\alpha)&=-2(\alpha-1)\frac{\langle \nabla^2 u_\alpha, \nabla u_\alpha \rangle \nabla u_\alpha}{1+|\nabla u_\alpha|^2}\nonumber \\
&=: f_\alpha. \label{SUeq2}
\end{align}
Now we let $h=h(r)$ be a piecewise linear function which equals the mean value of $u_\alpha$ on $\{\frac{R_2}{4}\} \times S^1$ and $\{2^k R_1\} \times S^1$ for all $1\le k\le l-1$. With the help of this we have
\begin{align*}
\Delta (u_\alpha-h) +A(u_\alpha)(\nabla u_\alpha, \nabla u_\alpha)=f_\alpha.
\end{align*}
Testing this equation with $u_\alpha-h$ and integrating over $A_k$ we get
\begin{align*}
\int_{A_k}|\nabla (u_\alpha-h)|^2=&\int_{A_k}(u_\alpha-h)(A(u_\alpha)(\nabla u_\alpha, \nabla u_\alpha)-f_\alpha)\\
&+2^{k+1}R_1\int_0^{2\pi}(u_\alpha-h)(u_\alpha-h)_r(2^{k+1}R_1,\theta)d\theta\\ &-2^kR_1\int_0^{2\pi}(u_\alpha-h)(u_\alpha-h)_r(2^kR_1,\theta)d\theta.
\end{align*}
We remark that the boundary integrals of $(u_\alpha-h)h_r$ vanish since $h$ is equal to the mean value of $u_\alpha$ on these boundaries and $h_r$ is piecewise constant. Because of \eqref{SUtangential2} and the Sobolev embedding (which we only apply on the annuli $A_k$) we know that for every $x\in A_k$ we have 
\begin{align}
|u_\alpha-h|(x)+\sum_{i=1}^2 |x|^i|\nabla^i u_\alpha| \le c\delta^{\frac{1}{2}}.\label{oscSU}
\end{align}
This implies that 
\begin{align*}
\int_{A_k}|\nabla (u_\alpha-h)|^2 \le& c \delta^{\frac{1}{2}}\int_{A_k}(|\nabla u_\alpha|^2+|f_\alpha|)\\
&+2^{k+1}R_1\int_0^{2\pi}(u_\alpha-h)(u_\alpha)_r(2^{k+1}R_1,\theta)d\theta\\ 
&-2^kR_1\int_0^{2\pi}(u_\alpha-h)(u_\alpha)_r(2^kR_1,\theta)d\theta.
\end{align*}
Taking the sum over $k$ we get
\begin{align*}
\int_{A(2R_1,\frac{R_2}{4})}|\nabla (u_\alpha-h)|^2 \le& c \delta^{\frac{1}{2}}\int_{A(2R_1,\frac{R_2}{4})}(|\nabla u_\alpha|^2+|f_\alpha|)\\
&+\frac{R_2}{4}\int_0^{2\pi}(u_\alpha-h)(u_\alpha)_r(\frac{R_2}{4},\theta)d\theta\\ 
&-2R_1\int_0^{2\pi}(u_\alpha-h)(u_\alpha)_r(2R_1,\theta)d\theta\\
\le& c\delta^{\frac{1}{2}}(1+(\log R_1^{1-\alpha})),
\end{align*}
where we used \eqref{oscSU} to estimate
\begin{align*}
\int_{A(2R_1,\frac{R_2}{4})}|f_\alpha| \le& c (\alpha-1) \int_{A(2R_1,\frac{R_2}{4})}|\nabla^2 u_\alpha|\\
\le& (\log R_1^{1-\alpha}).
\end{align*}
This finishes the proof of the Lemma.
\end{proof}
\subsection{Proof of the energy identity}
\begin{proof}
Because of the induction argument of Ding \& Tian \cite{ding95} we know that it is enough to prove the energy identity in the presence of one bubble. Since we are dealing with a local problem we assume from now on that $u_\alpha:\R^2 \supset B_1\rightarrow N$ and that we have only one energy concentration point $x^1=0$. Using the notations from section $2.1$ we assume that we obtain the bubble by rescaling with the factor $r_k^1=r_k$. From the smooth convergence $u_{\alpha_k}\rightarrow u_1$ away from $0$ we conclude that 
\begin{align*}
E_{\alpha_k}(u_{\alpha_k},B_1\backslash B_{R_0})\rightarrow E(u_1,B_1 \backslash B_{R_0})+\on{vol}(B_1\backslash B_{R_0}),
\end{align*}
for every $0<R_0<1$. Similarly, from the local $C^1$-convergence $v_k^1=v_k=u_{\alpha_k}(r_k\cdot)\rightarrow \omega$, we have for every $R>0$
\begin{align*}
E_{\alpha_k}(u_{\alpha_k}, B_{Rr_k})\rightarrow E(\omega). 
\end{align*}
Moreover this also implies that for every $R>0$ and $M>0$
\begin{align}
E_{\alpha_k}(u_{\alpha_k};B_{R_0}\backslash B_{\frac{R_0}{M}})+E_{\alpha_k}(u_{\alpha_k};B_{M r_k R}\backslash B_{r_k R}) &\to 0,\label{rest} 
\end{align}
as $k\rightarrow \infty$ and $R_0 \rightarrow 0$. Therefore it is easy to see that the proof of the energy identity in the case of one bubble is reduced to showing that
\begin{align}
\lim_{R\rightarrow \infty} \lim_{R_0 \rightarrow 0} \lim_{k\rightarrow \infty}  E_{\alpha_k}(u_{\alpha_k},A(Rr_k,R_0))= 0. \label{conv3}
\end{align}
Next we claim that due to the fact that we have only one bubble we can assume that for any $\delta>0$ there exists $k_0>0$ such that for all $k>k_0$ we have
\begin{align}
E(u_{\alpha_k},B_{2r}\backslash B_r) <\delta, \label{conv4}
\end{align}
for every $Rr_k\le r\le \frac{R_0}{2}$. To see this we argue by contradiction. If the claim is false, we may assume that as $k\to \infty$ there exists $s_k\in (Rr_k,\frac{R_0}{2})$ such that
\begin{align*}
E(u_{\alpha_k};B_{2s_k}\backslash B_{s_k})&=\max_{r\in (Rr_k,\frac{R_0}{2})}E(u_{\alpha_k};B_{2r}\backslash B_{r}) \nonumber \\
&\geq \delta.
\end{align*}
From \eqref{rest} we get that
\begin{align}
\frac{R_0}{s_k} &\to \infty \ \ \on{and} \nonumber \\
\frac{Rr_k}{s_k} &\to 0. \label{rk}
\end{align}
By defining 
\begin{align*}
\tilde{v}_k&:B_{\frac{R_0}{s_k}}\backslash B_{\frac{r_k R}{s_k}} \to N \nonumber \\
\tilde{v}_k(x)&=u_{\alpha_k}(s_k x) 
\end{align*}
we have that $\tilde{v}_k$ solves \eqref{SUeq} with $1$ replaced by $(s_k)^2$ and 
\begin{align}
\int_{B_{\frac{R_0}{s_k}}\backslash B_{\frac{r_k R}{s_k}}} ((s_k)^2+|\nabla \tilde{v}_k|^2)^{\alpha_k}&\leq cs_k^{2(\alpha_k-1)}, \label{bound} \\
E(\tilde{v}_k;B_2 \backslash B_1) &\geq \delta. \label{nontrivial}
\end{align}
By \eqref{bound}, \eqref{rk}, Proposition \ref{concent} and the arguments of section $2.1$ we may assume that $\tilde{v}_k \rightharpoonup \tilde{v}_0$ weakly in $W^{1,2}_{\on{loc}}(\R^2\backslash \{0 \},N)$, where $\tilde{v}_0:\R^2  \to N$ is a harmonic map with finite Dirichlet energy.
\newline
We have two possibilities. The first one is that there exists $\tilde{r}>0$ such that 
\begin{align*}
\sup_{k\in \N} \sup_{x\in B_{4} \backslash B_{\frac{1}{4}}} E(\tilde{v}_k;B_{\tilde{r}}(x)) < \varepsilon_0.
\end{align*}
With the help of Theorem \ref{regSU} and a covering argument this implies that $\tilde{v}_k\to \tilde{v}_0$ in 
\newline
$C^\infty(B_2 \backslash B_1,N)$. Since $\R^2 \backslash \{0\}$ is conformally equivalent to $S^2 \backslash \{N,S\}$ we conclude from \eqref{nontrivial} and the point removability result of Sacks \& Uhlenbeck \cite{sacks81}, that $\tilde{v}_0$ can be lifted to a smooth non-trivial harmonic map from $S^2$ to $N$, contradicting the assumption that we have only one bubble $\omega$.
\newline
The second possibility is that we have at least one energy-concentration point $y \in B_{4}\backslash B_{\frac{1}{4}}$. Now we can apply the blow-up procedure of section $2.1$ to conclude that there must exist a non-trivial harmonic two-sphere, again contradicting the assumption that there is only one bubble. This proves \eqref{conv4} and hence we can combine Theorem \ref{regSU}, Proposition \ref{concent}, Lemma \ref{radialSU} and Lemma \ref{SUtangential} (with $R_1=Rr_k$ and $R_2=R_0$) to estimate 
\begin{align*}
\int_{A(2Rr_k,\frac{R_0}{4})} (1+|\nabla u_{\alpha_k}|^2)^\alpha \le& c \int_{A(2Rr_k,\frac{R_0}{4})}|\nabla u_{\alpha_k}|^2+o_{R_0}(1)\\
=& c\int_{2Rr_k}^{\frac{R_0}{4}}\int_0^{2\pi} (r|(u_{\alpha_k})_r|^2+\frac{1}{r}|(u_{\alpha_k})_\theta|^2)dr d\theta+o_{R_0}(1)\\
\le& c\int_{2Rr_k}^{\frac{R_0}{4}}\int_0^{2\pi} \frac{1}{r}|(u_{\alpha_k})_\theta|^2dr d\theta+o_{R_0}(1)\\
&+ c(\alpha_k-1)\int_{2Rr_k}^{\frac{R_0}{4}}\frac{1}{r}(\int_{B_r}|\nabla u_{\alpha_k}|^{2\alpha_k} dx)dr\\
\le& o_k(1)+o_{R_0}(1)+c\sqrt{\delta}+c(1-\alpha_k)\log(Rr_k)\\
\le& o_k(1)+o_{R_0}(1)+c\sqrt{\delta} ,
\end{align*}
which, combined with \eqref{rest}, proves \eqref{conv3} (since $\delta>0$ was arbitrary) and therefore the main Theorem in the case of one bubble.
\end{proof}
\begin{remark}\label{optimal}
By a careful inspection of the above proof it is easy to see that the energy identity remains true for general sequences of critical points of $E_\alpha$ if and only if 
\[
\lim_{k\rightarrow \infty} (r_k^{i,j})^{1-\alpha_k}=1,
\]
for all $1\le i \le l$ and all $1\le j\le j_i$. This fact has also been observed by Li \& Wang \cite{li07}
\end{remark}

\section{Energy identity for the biharmonic approximation of harmonic maps}

In this section we prove Theorem \ref{main} for the biharmonic approximation of harmonic maps.
\subsection{Estimates and blow-up}
In the following we consider sequences of critical points $u_\eps\in C^\infty (M,N)$ ($\eps \rightarrow 0$) of the functional $E_\eps$ with uniformly bounded energy $E_\eps(u_\eps) \le c$ and which satisfy \eqref{epsent}. First of all we choose a subsequence $\eps_k \rightarrow 0$ such that 
\begin{align}
\eps_k \log(\frac{1}{\eps_k})\int_M |\Delta u_{\eps_k}|^2 =o_k(1).\label{epsent1}
\end{align}
Due to the uniform bound on the $W^{1,2}$-norm of $u_{\eps_k}$ we get the existence of a further subsequence (still denoted by $\eps_k$) such that $u_{\eps_k} \rightharpoonup u_0$ weakly in $W^{1,2}(M,N)$. In \cite{lamm06} we were able to show the following small energy estimate (see Corollary $2.10$ in \cite{lamm06}).
\begin{theorem}\label{sefourth}
There exists $\delta_0>0$ and $c>0$ such that if $u_\varepsilon \in C^\infty (M,N)$ is a solution of \eqref{fourth2} with $\int_{B_{2R}} (|\nabla u|^2+\varepsilon |\Delta u|^2)<\delta_0$ then we have for $\varepsilon$ small enough and every $m\in \N$
\begin{align}
\on{osc}_{B_R} u_\varepsilon +R^m||\nabla^m u_\varepsilon||_{L^\infty(B_R)} \le c( \int_{B_{2R}} (|\nabla u|^2+\varepsilon |\Delta u|^2))^{\frac{1}{2}}. \label{sefourth1}
\end{align}
\end{theorem}
Hence, as in section $2.1$, the sequence $u_\eps$ converges strongly to $u_0$ away from finitely many singular points $x^i\in M$, $1\le i\le l$, which are characterized by the condition
\begin{align}
\limsup_{k\rightarrow \infty} E_{\eps_k} (u_{\eps_k}, B_R(x^i)) \ge \delta_0, \label{econcbih}
\end{align}
for every $R>0$ and every $1\le i \le l$. Around these finitely many singular points we were able to perform a blow-up similar to the one of section $2.1$ (see section $3$ of \cite{lamm06}). Namely, for $R_0>0$ such that $B_{R_0}(x^i)\cap B_{R_0}(x^j) =\emptyset$ for every $1\le i\not= j \le l$, there exists a sequence of points $x^i_k\rightarrow x^i$ and a sequence of radii $r_k^i\rightarrow 0$ such that
\begin{align}
\max_{y\in B_{R_0}(x^i)}E_{\eps_k}(u_{\eps_k},B_{r_k^i}(y))=E_{\eps_k}(u_{\eps_k},B_{r_k^i}(x^i_k))=\frac{\delta_0}{2}. \label{fixradius}
\end{align}
Defining 
\begin{align}
 w^i_k:&B_{\frac{R_0}{r_k^i}}\rightarrow N, \nonumber \\
 w^i_k(x) &=u_{\eps_k}(x^i_k+r_k^i x) \label{defv}
\end{align}
we see that $w^i_k$ solves \eqref{fourth2} with $\eps_k$ replaced by $\tilde{\eps}_k=\frac{\eps_k}{(r_k^i)^2}$ and 
\begin{align}
\max_{y\in B_{\frac{R_0}{2r_k^i}}}E_{\tilde{\eps}_k}(w^i_k,B_1(y))=E_{\tilde{\eps}_k}(w^i_k,B_1(0))=\frac{\delta_0}{2}. \label{fixradius1}
\end{align}
Hence we can apply Theorem \ref{sefourth} to $w^i_k$ and conclude that $w^i_k$ converges smoothly to some map $\omega^i\in C^\infty \cap W^{1,2}(\R^2,N)$. Then we were able to show (Lemma $3.1$ in \cite{lamm06}) that for every $1\le i \le l$
\begin{align}
\tilde{\eps}_k=\frac{\eps_k}{(r_k^i)^2} \rightarrow 0, \label{concbih}
\end{align}
and therefore $\omega^i$ is a harmonic map with finite Dirichlet energy and can therefore be lifted to a smooth harmonic map from $S^2$ to $N$.

\subsection{Stress-energy tensor}
For a smooth map $u$ we have the well-known stress-energy tensor $S^1_{\alpha \beta}(u)$ given by
\begin{align}
S^1_{\alpha \beta}(u)=\frac{1}{2} |\nabla u|^2 \delta_{\alpha \beta}-\langle \nabla_\alpha u,\nabla_\beta u\rangle. \label{se1}
\end{align}
An easy calculation shows that if $u$ is a harmonic map then we have  
\begin{align}
\partial_\alpha S^1_{\alpha \beta}(u) =-\langle \Delta u,\nabla_\beta u\rangle=0. \label{se2}
\end{align}
Again for a smooth map $u$ we have the stress-energy tensor $S^2_{\alpha \beta}(v)$ defined by (see \cite{jiang87} and \cite{loubeau07})
\begin{align}
S^2_{\alpha \beta}(u)=\frac{1}{2} |\Delta u|^2 \delta_{\alpha \beta} +\langle \nabla_\gamma u, \nabla_\gamma \Delta u\rangle\delta_{\alpha \beta} -\langle \nabla_\alpha u, \nabla_\beta \Delta u\rangle-\langle \nabla_\beta u, \nabla_\alpha \Delta u\rangle. \label{se3}
\end{align}
By another easy calculation we see that if $u$ is an extrinsic biharmonic map (i.e. a solution of $\Delta^2 u \perp T_u N$) then we have
\begin{align}
\partial_\alpha S^2_{\alpha \beta} (u)= -\langle \nabla_\beta u, \Delta^2 u \rangle =0 . \label{se4}
\end{align}
Combining \eqref{se2} and \eqref{se4} we see that
\begin{align}
\partial_\alpha (S^1_{\alpha \beta}(u_\varepsilon )-\varepsilon S^2_{\alpha \beta}(u_\varepsilon) )=\langle \nabla_\beta u_\varepsilon, (\varepsilon \Delta^2-\Delta)u_\varepsilon \rangle =0, \label{se5}
\end{align}
if $u_\varepsilon$ is a solution of \eqref{fourth2}. As in the case of harmonic maps (see \cite{sacks81}) we use this divergence-free quantity to get a Hopf differential type estimate for solutions of \eqref{fourth2}.
\begin{lemma}\label{hopffourth}
Let $u_\varepsilon \in C^\infty(B,N)$ be a solution of \eqref{fourth2}. Then we have for all $0<r<1$
\begin{align}
\int_{\partial B_r}|(u_\eps)_r|^2 \le \frac{1}{r^2}\int_{\partial B_r} |(u_\eps)_\theta|^2 +\frac{c\varepsilon}{r} \int_{B_r} |\Delta u_\eps|^2 +c\varepsilon \int_{\partial B_r} (|\Delta u_\eps|^2 +|\nabla u_\eps||\nabla^3 u_\eps|). \label{hopffourth1}
\end{align}
\end{lemma}
\begin{proof}
Multiplying \eqref{se5} by $x^\beta$ and integrating by parts we get for every $0<r<1$
\begin{align*}
 \int_{B_r} (S^1_{\alpha \beta}(u_\varepsilon)-\varepsilon S^2_{\alpha \beta}(u_\varepsilon )) \delta_{\alpha \beta} =\int_{\partial B_r} (S^1_{\alpha \beta}(u_\varepsilon)-\varepsilon S^2_{\alpha \beta}(u_\varepsilon))x^\beta \nu ^\alpha,
\end{align*}
where $\nu$ is the outer unit normal to $\partial B_r$. Now we calculate
\begin{align*}
(S^1_{\alpha \beta}(u_\varepsilon)-\varepsilon S^2_{\alpha \beta}(u_\varepsilon))\delta_{\alpha \beta} =-\varepsilon |\Delta u_\varepsilon|^2
\end{align*}
and
\begin{align*}
(S^1_{\alpha \beta}(u_\varepsilon)-\varepsilon S^2_{\alpha \beta}(u_\varepsilon))x^\beta \nu^\alpha =& \frac{r}{2} |\nabla u_\varepsilon|^2-r|(u_\varepsilon)_r|^2\\
&-r\varepsilon(\frac{1}{2}|\Delta u_\varepsilon|^2 +\langle \nabla u_\varepsilon,\nabla \Delta u_\varepsilon \rangle  -2 \langle (u_\varepsilon)_r,(\Delta u_\varepsilon)_r \rangle )\\
=&\frac{1}{2r}|(u_\varepsilon)_\theta|^2 -\frac{r}{2} |(u_\varepsilon)_r|^2\\
&-r\varepsilon(\frac{1}{2}|\Delta u_\varepsilon|^2 +\langle \nabla u_\varepsilon,\nabla \Delta u_\varepsilon \rangle  -2 \langle (u_\varepsilon)_r,(\Delta u_\varepsilon)_r \rangle ),
\end{align*}
where we used the identity $|\nabla u|^2 =|u_r|^2+\frac{1}{r^2}|u_\theta|^2$. This finishes the proof of the Lemma.
\end{proof}
\subsection{Estimate for the tangential component}
In this subsection we prove an estimate for the biharmonic approximation similar to the one given in section $2.3$ for the Sacks-Uhlenbeck approximation.
\begin{lemma}\label{fourthtangential}
There exists $\delta_1>0$ such that for all $\delta < \delta_1$ and all solutions $u_\varepsilon$ of \eqref{fourth2} with $\int_{A(r,2r)}(|\nabla u_\varepsilon|^2+\varepsilon |\Delta u_\varepsilon|^2) < \delta$ for every $r\in (R_1,\frac{R_2}{2})$, we have for $\varepsilon$ small enough
\begin{align}
\int_{2R_1}^{\frac{R_2}{4}}\int_0^{2\pi} \frac{1}{r} |(u_\varepsilon)_\theta|^2 dr d\theta \le c\sqrt{\delta}(1+\frac{\varepsilon}{(R_1)^2}). \label{fourthtangential1}
\end{align}
\end{lemma}
\begin{proof}
The proof follows directly from the one of Lemma \ref{SUtangential}. Namely instead of using Theorem \ref{regSU} we use Theorem \ref{sefourth} to conclude that
\begin{align}
\sum_{i=1}^4 |x|^i |\nabla^i u_\varepsilon| \le c\sqrt{\delta} \label{tangbi1}
\end{align}
for every $x\in A(2R_1,\frac{R_2}{4})$. Moreover we note that equation \eqref{fourth2} can equivalently be written as
\begin{align}
\Delta u_\varepsilon +A(u_\varepsilon)(\nabla u_\varepsilon, \nabla u_\varepsilon) &=\varepsilon \Delta^2 u_\varepsilon +f[u_\varepsilon]\nonumber \\
&= f_\varepsilon. \label{tangbi2}
\end{align}
Using this form of the equation it is easy to see that the proof of Lemma \ref{SUtangential} carries over to this situation once we notice that because of \eqref{fourth3} and \eqref{tangbi1} we have
\begin{align*}
\sqrt{\delta} \int_{A(2R_1,\frac{R_2}{4})}|f_\varepsilon| &\le c\sqrt{\delta} \varepsilon \int_{A(2R_1,\frac{R_2}{4})}(|\nabla^4 u_\varepsilon|+|\nabla u_\varepsilon||\nabla^3 u_\varepsilon|+|\nabla^2 u_\varepsilon|^2+|\nabla u_\varepsilon|^4)\\
&\le c\sqrt{\delta} \frac{\varepsilon}{(R_1)^2}.
\end{align*}  
\end{proof}
\subsection{Proof of the energy identity}
\begin{proof}
Following the remarks of section $2.4$ (using the results of section $3.1$) we can assume that we have only one energy concentration point $x^1=0 \in B_1 \subset \R^2$ and one bubble $\omega^1$ which is obtained by rescaling $u_{\varepsilon_k}$ by the factor $r_k^1=r_k$. Again the proof of the energy identity is reduced to showing that
\begin{align}
\lim_{R\rightarrow \infty} \lim_{R_0 \rightarrow 0} \lim_{k\rightarrow \infty} E_{\varepsilon_k}(u_{\eps_k} ,B_{R_0} \backslash B_{Rr_k}) =0. \label{convbih}
\end{align}
Using similar arguments as in section $2.4$ we can moreover assume that for any $\delta>0$ there exists $k_0>0$ such that for all $k>k_0$ we have
\begin{align}
E_{\eps_k}(u_{\eps_k},B_{2r}\backslash B_r) < \delta, \label{convbih2}
\end{align}
for every $Rr_k\le r\le \frac{R_0}{2}$. Hence we can apply \eqref{tangbi1} with $R_1=Rr_k$ and $R_2=R_0$ to get
\begin{align}
\eps_k \int_{A(2Rr_k,\frac{R_0}{4})} |\Delta u_{\eps_k}|^2 \le& c\delta \eps_k \int_{A(2Rr_k,\frac{R_0}{4})} \frac{dx}{|x|^4} \nonumber \\
\le& c\delta \frac{\eps_k}{R^2r_k^2} \nonumber \\
=& o_k(1), \label{convbih3}
\end{align}
where we used \eqref{concbih} in the last line. Combining \eqref{epsent1}, Lemma \ref{hopffourth}, Lemma \ref{fourthtangential}, \eqref{convbih2}, \eqref{convbih3} and \eqref{concbih} we get
\begin{align}
E_{\eps_k}(u_{\eps_k},A(2Rr_k,\frac{R_0}{4}))\le& \int_{Rr_k}^{\frac{R_0}{4}} \int_0^{2\pi} (r|(u_{\eps_k})_r|^2+\frac{1}{r}|(u_{\eps_k})_\theta|^2)dr d\theta +o_k(1)\nonumber \\
\le& c\int_{Rr_k}^{\frac{R_0}{4}} \int_0^{2\pi}\frac{1}{r}|(u_{\eps_k})_\theta|^2 dr d\theta +c\eps_k \int_{Rr_k}^{\frac{R_0}{4}}(\frac{1}{r}\int_{B_r} |\Delta u_{\eps_k}|^2)dr\nonumber \\
&+c\eps_k \int_{Rr_k}^{\frac{R_0}{4}}\int_{\partial B_r} (|\Delta u_{\eps_k}|^2+|\nabla u_{\eps_k}||\nabla^3 u_{\eps_k}|)+o_k(1)\nonumber \\
\le& c\eps_k \log(\frac{1}{Rr_k}) \int_M |\Delta u_{\eps_k}|^2+ \frac{c\eps_k}{(r_k)^2}+o_k(1)+c\sqrt{\delta}\nonumber \\
\le& c\eps_k \log(\frac{1}{\eps_k}) \int_M |\Delta u_{\eps_k}|^2+ o_k(1)+c\sqrt{\delta}\nonumber \\
\le&  o_k(1)+c\sqrt{\delta}, \label{convbih4}
\end{align}
which proves \eqref{convbih}.
\end{proof} 
\section{Proof of Lemma \ref{entropy}}

We follow closely the work of Struwe \cite{struwe00a} (see also \cite{struwe88a}, \cite{struwe88b} and \cite{struwe00}).  
\begin{proof}
Since the methods are very similar for both approximations we only prove the Lemma for $E_\alpha$.  
\newline
First of all we note that the minimax principle (see for example \cite{struwe00}, Theorem $4.2$) guarantees the existence of a critical point $u_\alpha$ of $E_\alpha$ with $E_\alpha(u_\alpha)=\beta_\alpha$. The difficult part now consists of showing that we can also find a sequence of critical points satisfying \eqref{ent}.
\newline
We note that it is easy to see that the function 
\begin{align*}
\alpha \rightarrow \beta_\alpha=\inf_{F\in \mathcal{F}}  \sup_{u\in F} E_\alpha(u)
\end{align*}
is non-decreasing and hence differentiable almost everywhere with differential $0\le \frac{d\beta_\alpha}{d\alpha} \in L^1([1,\alpha_1])$ for $\alpha_1>1$. Therefore it follows that
\begin{align}
2B=\liminf_{\alpha \rightarrow 1} (\alpha-1) \log(\frac{1}{\alpha-1}) \frac{d\beta_\alpha}{d\alpha}  =0. \label{entropy1}
\end{align}
To see this we assume that $B>0$ and we get for $(A-1)$ very small
\begin{align*}
\int_1^A \frac{d\beta_\alpha}{d\alpha} d\alpha \ge -B \int_1^A \frac{d\alpha}{(\alpha-1) \log(\alpha-1)}=\infty,
\end{align*}
which contradicts the fact that $\frac{d\beta_\alpha}{d\alpha} \in L^1([1,\alpha_1])$. Next we let $\alpha>1$ be a point of differentiabilty of $\beta_\alpha$ and we choose a sequence $\alpha_k \rightarrow \alpha$ ($\alpha_{k+1}\le \alpha_k$). For every $k\in \N$ we choose $F_k\in \mathcal{F}$ such that
\begin{align*}
\sup_{u\in F_k} E_{\alpha_k}(u) \le \beta_{\alpha_k} +(\alpha_k-\alpha). 
\end{align*}
Since $\beta_\alpha$ is differentiable in $\alpha$ we get that for sufficiently large $k$ we have
\begin{align*}
\beta_{\alpha_k} \le \beta_\alpha +(\frac{d\beta_\alpha}{d\alpha}+1)(\alpha_k-\alpha). 
\end{align*}
Combining the above two estimates we get
\begin{align*}
\beta_\alpha \le& \sup_{u\in F_k} E_\alpha (u) \le \sup_{u\in F_k} E_{\alpha_k}(u) \le \beta_{\alpha_k}+(\alpha_k-\alpha) \le \beta_\alpha +(\frac{d\beta_\alpha}{d\alpha}+2)(\alpha_k-\alpha). 
\end{align*}
Next we choose $v\in F_k$ such that
\begin{align*}
\beta_\alpha -(\alpha_k-\alpha) \le E_\alpha (v). 
\end{align*}
Combining all this gives the existence of a map $v$ such that
\begin{align}
\beta_\alpha-(\alpha_k-\alpha)\le& E_\alpha (v) 
\le E_{\alpha_k}(v)
\le \sup_{u\in F_k} E_{\alpha_k} (u)  \le \beta_{\alpha_k}+\alpha_k-\alpha  \nonumber \\
\le& \beta_\alpha+(\frac{d\beta_\alpha}{d\alpha}+2)(\alpha_k-\alpha). \label{entropy6}
\end{align}
Now we prove three intermediate steps.
\newline
\underline{Step $1$:} For every $v\in W^{1,2\alpha_k}(M,N)$ which satisfies \eqref{entropy6} we have the estimate
\begin{align}
\partial_\alpha E_\alpha(v) \le \frac{d\beta_\alpha}{d\alpha}+3.\label{entro1}
\end{align}
From \eqref{entropy6} we get
\begin{align*}
\frac{E_{\alpha_k}(v)-E_\alpha(v)}{\alpha_k-\alpha} \le \frac{d\beta_\alpha}{d\alpha}+3 
\end{align*}
and hence by the mean value theorem there exists a number $\alpha \le \alpha`\le \alpha_k$ such that
\begin{align}
\partial_\alpha E_{\alpha`}(v) \le \frac{d\beta_\alpha}{d\alpha}+3. \label{entropy8}
\end{align}
Since moreover 
\begin{align*}
\partial_\alpha E_\alpha(u)=&\int_M \log(1+|\nabla u|^2)(1+|\nabla u|^2)^\alpha \\
\le& \int_M \log(1+|\nabla u|^2)(1+|\nabla u|^2)^{\alpha`}\\
=& \partial_\alpha E_{\alpha`}(u)
\end{align*}
for every $u\in W^{1,2\alpha_k}(M,N)$ we finish the proof of step $1$.
\newline
\underline{Step $2$:} We have
\begin{align}
\sup \{ |\langle dE_{\alpha_k}(u),v \rangle -\langle dE_\alpha(u),v\rangle| ; ||v||_{W^{1,2\alpha_k}(u^\star TN)}\le 1 \} \rightarrow 0, \label{entro2}
\end{align}
where
\begin{align*}
W^{1,2\alpha_k}(u^\star TN)=\{v\in W^{1,2\alpha_k}(M,R^m)| v(x)\in T_{u(x)}N\ \ \forall \ \ x\in M\},
\end{align*} 
uniformly for all $u\in W^{1,2\alpha_k}(M,N)$ satisfying \eqref{entropy6}.
\newline
To see this we note that for every $v\in W^{1,2\alpha_k}(u^\star TN)$ with $||v||_{W^{1,2\alpha_k}(u^\star TN)}\le 1$ we have
\begin{align*}
|\langle &dE_{\alpha_k}(u),v \rangle -\langle dE_\alpha(u),v\rangle| \\
\le& \int_M (2\alpha_k(1+|\nabla u|^2)^{\alpha_k-1}-2\alpha(1+|\nabla u|^2)^{\alpha-1})|\nabla u||\nabla v|\\
=:& I.
\end{align*}
Now we estimate 
\begin{align*}
I \le& 2(\alpha_k-\alpha)(\int_M (1+|\nabla u|^2)^{\alpha_k})^{\frac{\alpha_k-1}{\alpha_k}}(\int_M |\nabla u|^{2\alpha_k})^{\frac{1}{2\alpha_k}} (\int_M |\nabla v|^{2\alpha_k})^{\frac{1}{2\alpha_k}}\\
&+2\alpha \int_M  ((1+|\nabla u|^2)^{\alpha_k-1}-(1+|\nabla u|^2)^{\alpha-1})|\nabla u||\nabla v|\\
&\le c(\alpha_k-\alpha) +2\alpha \int_M  ((1+|\nabla u|^2)^{\alpha_k-1}-(1+|\nabla u|^2)^{\alpha-1})|\nabla u||\nabla v|.
\end{align*}
Next we use the estimate $2|\nabla u||\nabla v| \le \frac{1}{\delta} (1+|\nabla u|^2)+\delta |\nabla v|^2$, \eqref{entropy6} and Young's inequality to get
\begin{align*}
2\alpha \int_M  ((1+|\nabla u|^2)^{\alpha_k-1}&-(1+|\nabla u|^2)^{\alpha-1})|\nabla u||\nabla v|\\
\le& \frac{\alpha}{\delta} (E_{\alpha_k}(u)-E_\alpha (u))+\alpha \delta \int_M (1+|\nabla u|^2)^{\alpha_k-1}|\nabla v|^2\\
\le& \frac{c(\alpha_k-\alpha)}{\delta} +\delta \alpha(\frac{\alpha_k-1}{\alpha_k}E_{\alpha_k}(u)+\frac{1}{\alpha_k}\int_M |\nabla v|^{2\alpha_k}).
\end{align*}
Choosing $\delta=\sqrt{\alpha_k-\alpha}$ we conclude that
\begin{align*}
I\le c\sqrt{\alpha_k-\alpha} \rightarrow 0
\end{align*}
and this proves \eqref{entro2}.
\newline
\underline{Step $3$:} There exists a sequence $u_k\in W^{1,2\alpha_k}(M,N)$ satisfying \eqref{entropy6} and
\begin{align}
||dE_{\alpha_k} (u_k)||_{(W^{1,2\alpha_k}(M,N))^\star} \rightarrow 0. \label{entropy9}
\end{align}
If this is not the case we can find $\delta>0$ such that
\begin{align*}
||dE_{\alpha_k} (u)||_{(W^{1,2\alpha_k}(M,N))^\star} \ge 4\delta 
\end{align*}
for all $u$ satisfying \eqref{entropy6} and all $k$ large enough. For these $k$ we let 
\newline
$e_k:W^{1,2\alpha_k}(M,N) \rightarrow W^{1,2\alpha_k}(u^\star TN)$ be a locally Lipschitz continuous pseudo-gradient vectorfield for $E_{\alpha_k}$ with $||e_k(u)||_{W^{1,2\alpha_k}(u^\star TN)} \le 1$ and
\begin{align*}
\langle dE_{\alpha_k}(u) ,e_k(u)\rangle \le -\frac{1}{2} ||dE_{\alpha_k}(u)||_{(W^{1,2\alpha_k}(M,N))^\star}\le -2\delta,
\end{align*}
for all $u$ satisfying \eqref{entropy6}.
\newline
Let $\psi \in C^\infty(\R)$ be cut-off function such that $0\le \psi \le 1$, $\psi(s)=0$ for $s\le 0$, $\psi(s)=1$ for $s\ge 1$ and for $k$ large enough we let
\begin{align*}
\psi_k(u)=\psi \Big(\frac{E_\alpha(u) -\big(\beta_\alpha-(\alpha_k-\alpha)\big)}{\alpha_k-\alpha}\Big).
\end{align*}
Since $e_k$ is Lipschitz continuous the vectorfield
\begin{align*}
\tilde{e}_k(u)=\psi_k(u) e_k(u)
\end{align*}
then also defines a Lipschitz continuous tangent vectorfield. Finally we let $\phi_k:\R_0^+\times W^{1,2\alpha_k}(M,N)\rightarrow  W^{1,2\alpha_k}(M,N)$ be the flow generated be $\tilde{e}_k$:
\begin{align}
\frac{d}{dt} \phi_k(t,u) &=\tilde{e}_k(\phi_k(t,u)),\ \ \ t>0\nonumber \\
\phi_k(0,u)&=u. \label{entropy12}
\end{align}
Let $F_k\in \mathcal{F}$ be chosen as above and define for $v \in F_k$,  $v_t=\phi_k(t,v)$. Then we know from the assumptions of the Lemma that $v_t \in F_k$ for all $t\in \R^+_0$ and that
\begin{align*}
\sup_{v\in F_k} E_{\alpha_k}(v_t)\le  \sup_{v\in F_k} E_{\alpha_k}(v)\le \beta_{\alpha_k}+(\alpha_k-\alpha) 
\end{align*}
for all $t\ge 0$. Hence 
\begin{align}
M(t)=\sup_{v\in F_k} E_\alpha (v_t)\ge \beta_\alpha \label{entropy15}
\end{align}
is attained only at points $v_0$ for which $(v_0)_t$ satisfies \eqref{entropy6}. By noting that this implies $\psi_k((v_0)_t)=1$ we calculate 
\begin{align*}
\frac{d}{dt} E_{\alpha}((v_0)_t) =& \langle dE_{\alpha}((v_0)_t), \frac{d}{dt} (v_0)_t \rangle \\
=& \psi_k((v_0)_t)\langle dE_{\alpha}((v_0)_t),e_k((v_0)_t)\rangle \\
\le& \langle dE_{\alpha_k}((v_0)_t),e_k((v_0)_t)\rangle+|\langle dE_\alpha((v_0)_t)-dE_{\alpha_k}((v_0)_t),e_k((v_0)_t) \rangle|\\
\le& -2\delta +o_k(1), 
\end{align*}
where we used \eqref{entro2} in the last step. This shows that for $k$ large enough we get
\begin{align}
\frac{d}{dt} M(t) \le -\delta <0
\end{align}
and hence $M(t) <\beta_\alpha$ for large $t$ contradicting the definition of $\beta_\alpha$. Altogether this finishes the proof of step $3$.
\newline
To finish the proof of the Lemma we consider a sequence $u_k\in W^{1,2\alpha_k}(M,N)$ satisfying \eqref{entropy6} and \eqref{entropy9}. We know that $||u_k||_{W^{1,2\alpha_k}(M,N)}\le c$ and therefore we may assume that $u_k \rightharpoonup u_\alpha$ weakly in $W^{1,2\alpha}(M,N)$ and strongly in $L^{2\alpha}\cap C^{0,\beta}(M,N)$ for some $0<\beta <1$. Since $C^\infty(M,N)$ is dense in $W^{1,2\alpha}(M,N)$ we can moreover find a sequence $u^l\in C^\infty(M,N)$ such that $u^l\rightarrow u_\alpha$ strongly in $W^{1,2\alpha}(M,N)$. 
\newline
Next we define the functional $F_\alpha: W^{1,2\alpha}(M,\R^m) \rightarrow \R$ by
\begin{align}
F_\alpha (u) =\int_M (1+|\nabla u|^2)^\alpha. \label{modenergy}
\end{align}
Clearly we have that 
\begin{align*}
E_\alpha(u)=F_\alpha(u)
\end{align*}
for all $u\in W^{1,2\alpha}(M,N)$. For $v\in W^{1,2\alpha}(M,N)$ we define the projection 
\begin{align*} 
P_v:W^{1,2\alpha}(M,\R^m)&\rightarrow W^{1,2\alpha} (v^\star TN).
\end{align*}
Following the proof of Lemma $3.26$ in \cite{urakawa93} we get that
\begin{align*}
||(id -P_{u_k})(u_k-u^l)||_{W^{1,2\alpha_k}(M,\R^m)}\rightarrow 0,
\end{align*}
as $k,l\rightarrow \infty$. Hence we get from \eqref{entropy9} as in Lemma $3.7$ of \cite{urakawa93} that 
\begin{align}
| \langle d F_{\alpha_k} (u_k), u_k-u^l \rangle | \rightarrow 0,\label{convi1}
\end{align}
as $k,l\rightarrow \infty$. By convexity we know that for every $k,l\in \N$ we have
\begin{align*}
\int_M (1+|\nabla u^l|^2)^{\alpha_k}\ge& \int_M (1+|\nabla u_k|^2)^{\alpha_k}+\alpha_k \int_M (1+|\nabla u_k|^2)^{\alpha_k-1}\nabla u_k \nabla(u^l-u_k)\\
&+\alpha_k \int_M|\nabla(u_k-u^l)|^2.
\end{align*}
For any fixed $l\in \N$ we use this together with \eqref{convi1} to get
\begin{align*}
o_k(1)=& \langle dF_{\alpha_k}(u_k),u_k-u^l \rangle\\
\ge&  \int_M \Big( (1+|\nabla u_k|^2)^{\alpha_k} -(1+|\nabla u^l|^2)^{\alpha_k}+\alpha_k|\nabla(u_k-u^l)|^2 \Big)\\
\ge& \int_M \Big( (1+|\nabla u_k|^2)^{\alpha} -(1+|\nabla u^l|^2)^{\alpha}+\alpha|\nabla(u_k-u^l)|^2 \Big)\\
&- c(||\nabla u^l||_{L^\infty(M,\R^m)})(\alpha_k-\alpha).
\end{align*}
By letting first $k\rightarrow \infty$ and then $l\rightarrow \infty$ we conclude that $\nabla u_k \rightarrow \nabla u_\alpha$ pointwise a.e. and $E_\alpha(u_k)\rightarrow E_\alpha (u_\alpha)$. Hence we finally get that $u_k \rightarrow u_\alpha$ strongly in $W^{1,2\alpha}(M,N)$. By \eqref{entropy6} we have $E_{\alpha_k}(u_k)\rightarrow E_\alpha(u_\alpha)=\beta_\alpha$ and by step $2$ and $3$ we conclude that $u_\alpha$ is a critical point of $E_\alpha$. Since the function $s  \rightarrow \log(1+s^2)(1+s^2)^\alpha$ is convex we know that $\partial_\alpha E_\alpha$ is lower semi-continuous on $W^{1,2\alpha}(M,N)$ and therefore we can use step $1$ to get
\begin{align*}
\partial_\alpha E_\alpha(u_\alpha)\le \liminf_{k\rightarrow \infty} \partial_\alpha E_\alpha (u_k)\le \frac{d\beta_\alpha}{d\alpha}+3.
\end{align*}
Combining all this with \eqref{entropy1} we finish the proof of the Lemma.
\end{proof}

\bibliographystyle{amsplain}

\end{document}